\documentclass[12pt,a4paper,twoside]{article}
\usepackage{graphicx}
\usepackage{latexsym}
\usepackage{amssymb}
\newtheorem{thm}{Theorem}
\newcommand{\E}{I\kern-0.37emE}
\newcommand{\Q}{I\kern-0.37emP}

\newcommand{\ny}{n\rightarrow\infty}
\newcommand{\boldgreek}[1]{\mbox{\boldmath$#1$}}

\begin{document}

\title{Rank Tests for Corrupted Linear Models\thanks{P. K. Sen was partially supported by the Boshamer Professorship Research
funding at the University of North Carolina, Chapel Hill. The work of J. Jure\v{c}kov\'a was partially 
supported by the Grant GA\v{C}R 201/12/0083 and the work of Jan Picek was partially supported by the Grant GA\v{C}R P209/10/2045. }
}


\author{Pranab Kumar Sen$^{\rm a,\,}$\thanks{Corresponding author. E-mail: {\tt pksen@bios.unc.edu}
}       \and
        Jana Jure\v{c}kov\'a$^{\rm b}$ \and 
        Jan Picek$^{\rm c}$\\[5mm]
$^{\rm a}${\em{University of North Carolina, Chapel Hill, USA}}
\\$^{\rm b}${\em Charles University in Prague, Czech Republic}
\\$^{\rm c}${\em Technical University of Liberec, Czech Republic}
 }



\date{\small{\textit{Journal of the Indian Statistical Association} \textit{51/1, 2013, 201-229}}  }

\maketitle

\begin{abstract}
For some variants of regression models, including partial, measurement error or error-in-variables, latent effects, semi-parametric and otherwise corrupted linear models, the classical parametric tests generally do not perform well. Various modifications and generalizations considered extensively in the literature rests on stringent regularity assumptions which are not likely to be tenable in many applications. However, in such non-standard cases, rank based tests can be adapted better, and further, incorporation of rank analysis of covariance tools enhance their power-efficiency. Numerical studies and a real  data illustration show the superiority of rank based inference in such corrupted linear models.
\end{abstract}


\section{Introduction}
\label{intro}

Classical linear regression models induce some stringent additivity, linearity, homoscedasticity and normality assumptions which may not be tenable in many applications giving rise to the so called corrupted linear models where one or more of these assumptions may not be tenable. In simple nonparametric linear models, the normality assumption has been dispensed with in favor of a more general class of continuous distributions. Yet, in more contemporary applications in biomedical, clinical and genomics studies, the very assumption of linearity may be questionable. Sans such a linear setup, the performance of rank based testing procedures may be generally far better than their strict parametric counterparts. Our contemplated corrupted linear models relate to this scenario where the basic linearity assumption is vitiated by possible error-in-variables, measurement errors, possible latent effects, and the so called random effects and mixed effects; even partial linear models and some semi-parametric models belong to this contemplated class. For example, Fuller (1987) has detailed a large class of models which can be classified as measurement error or error-in-variable models; some genuine identifiability issues may crop-up in the use of standard parametric inference. Another variation is the usual regression models with stochastic predictors whose possible non-normal distribution can create stumbling blocks to the adaption of standard parametric methods. In addition such stochastic predictors may not be linearly related with the primary response variable. The impact of such nonregular setups on statistical tests has been 
considered by Ghosh and Sen (1971), followed by more general treatise by others. In a semi-parametric setup, partial linear models were introduced 
mostly during the 1980s and 1990s (Heckman (1986), Speckman (1988), Khuri, Mathew and Sinha (1988),  Chen (1988), Gao (1995), He and Shi (1996), Liang et al. (1999), Hardle et al. (2000), He and Liang (2000), and Boente and Rodriguez 2006, among others).  Incorporation of measurement errors in this setup evolved first in nonlinear models (Carroll et al.\,(2006)) and then  in nonparametric setups  only in the past decade. For nonlinear models one may try to mimic the linear model setups with linear or quadratic approximations, but again those may call for a second source of non-robustness arising from such possibly inadequate approximations. Motivated by this diversity of models and the need for a unified view of such nonstandard or corrupted linear models, the present study mainly aims to introduce such corrupted linear models in a more general setup, exhibit the supremacy of rank based tests and illustrate its adaptability in some real applications.

%

Consider a {\em semiparametric partially 
linear} model where a real response Y is regressed
to a set of observable covariates ${\mathbf x}$ and further depends on some possibly 
unobservable $\mathbf Z$ in the form:
\begin{equation} \label{model}
 Y_{i}=\beta_0+{\bf x}_{i}^{\top}{\boldgreek\beta}+\nu(\mathbf Z_{i})+e_{i},
\quad i=1,\ldots,n,
\end{equation}
where the ${\mathbf x}_{i}$ are known (non-stochastic) $p$-vectors, not all the same,
$\mathbf{Z}_{i}$ is a stochastic $q$-vector covariate $(q \geq 1),$ and
the form of the function $\nu(\mathbf{Z}_{i})$ is unspecified. Moreover, the
$\mathbf Z_{i}$ may be observable, partially observable or unobservable;
in the latter case, they lead to latent effects models.
If the $\mathbf Z_{i}$ are observable, eventually with measurement errors, (\ref{model}) relates to a partially linear and measurement error
model. he unknown $\nu(.)$ links (\ref{model}) to the 
semiparametric model. A big advantage of the rank procedure is that it avoids a nonparametric estimation of unknown $\nu(.).$ 
The literature recommends the functional estimation procedures, using
various smoothing tools; but they demand smoothness assumptions, while they
usually result in slower rates of convergence than the rank procedures. We refer to Heckman (1986), Speckman (1988), Chen (1988), He and Shi (1996), He and Liang (2000),  
Bianco et al. (2006), Boente and Rodriguez (2006),  
among other works. The H\"{a}rdle et al. (2000) monograph is noteworthy in this context.

An alternative approach is a transformation of variables in regression problems which achieves linearity or normality; but this usually sacrifices the homoscedasticity condition. The heteroscedastic models and models with measurement errors were intensively treated in the literature; we refer to  
Fuller (1987), Cheng and van Ness (1999)
and Carroll et al. (2006), and to additional references cited therein.

In contrast to the above methods, we put the main emphasis on nonparametric tests based on 
 rank statistics. They are valid also for non-normal error distributions, do not demand the finite variances, and their asymptotic forms typically have the standard rate of convergence $n^{-1/2}.$ The smoothing techniques as B-splines and kernel smoothing, which are commonly used for estimation in the
semiparametric linear models, generally require a large $n$ and result in a slower rate of convergence than $n^{-1/2}.$ 

 The problem of testing
the monotonicity of regression was considered by \cite{Ghosal2000}, who used a nonparametric approach in a semiparametric setup. These models can be sometimes reduced to 
(\ref{model}) by suitable reformulation.

   We often want to test the null hypothesis of no or partial regression of $Y$ on $\mathbf x,$ treating $\beta_0$ and $\nu(\cdot)$ as nuisance parameters
and functions, respectively. The statistical interest is then confined to the fixed-effect parameter $\boldgreek\beta,$ regarding $\nu(\cdot)$ as a nuisance function, similarly as in the \cite{Cox1972} proportional hazard model. More precisely, we want to test 
\begin{equation}\label{H0}
\mathbf H_0:~\boldgreek\beta=\mathbf 0 \quad vs \quad \mathbf H_1:~\boldgreek\beta\neq\mathbf 0
\end{equation}
with nuisance $\beta_0$ and $\nu(\cdot).$

Although $\nu(\cdot)$ is unspecified in (\ref{model}), it is of interest to 
distinguish two cases according as the covariate $\mathbf Z$ is observable or 
not. If $\mathbf Z$ is unobservable, (\ref{model}) corresponds to the \textit{latent effects}
model, although in the usual linear model setup, $\nu(\mathbf Z)$ is taken to be a linear functional, 
whereas in (\ref{model}) it is unspecified.   If  $\mathbf Z_i$'s are observable and regarded as identically distributed
random variables with some unspecified distribution and 
independent of the error $e_i$, then letting $e_i^* = e_i + \nu(\mathbf Z_i)$ we may 
still claim that the $e_i^*$ are i.i.d. random variables. However, their distribution function
is unlikely to be normal even if the $\mathbf Z_i$ were normally distributed;
this is specially because of the unspecified nature of $\nu(\cdot).$ 
 On the other hand, since the $e_i^*$  are   independent  identically distributed
random variables,  the classical nonparametric rank based tests are adaptable.
This naturally suggests that nonparametric tests based on rank statistics would 
have better scope as well as power properties.

%
There is a much better perspective if the $\mathbf Z_i$, though stochastic, are observable. Unlike the parametric analysis of covariance (ANOCOVA) the assumption of linearity of regression is not necessary in the nonparametric ANOCOVA approach. Quade (1969) considered a rank ANOCOVA procedure based on the rank sum statistics, and that was  extended immediately to general scores tests in more general models by Puri and Sen (1971) where earlier references are also cited. Even the work of Ghosh and Sen (1971) is closely related to this aspect of rank tests. In this context, by virtue of the fact that ranks are invariant under any strictly monotone transformation on the covariates, the linearity of the regression on covariate may no longer be necessary, and the resulting rank ANOCOVA tests are therefore much more robust than their parametric counterparts and typically have greater power than nonparametric rank ANOVA tests which ignore the covariates. This improvement comes out of the fact that the joint distribution of the coordinatewise rank statistics is typically close to a multinormal one and that validates the use of ANOCOVA tools even when the underlying form of $\nu(.)$ is nonlinear. Even more, the rank tests are still applicable if the $\mathbf Z_i$ are observable, 
but subject to measurement errors as in the model considered by Nummi and M\"{o}tt\"{o}nen (2004);  
then the $e_i^*$ are
still i.i.d. random variables, though with some other distribution function. This shows an advantage of the nonparametric analysis of covariance procedures comparing with other methods.

In a general regression setup where the regressors are stochastic, Ghosh and Sen (1971) modified the usual rank tests for testing the hypothesis of no regression, and in the measurement error model, Jure\v{c}kov\'a et al. (2010) considered suitable rank tests. In both the cases, the hypothesis of no regression generates the same invariance structure which validates  the conventional rank tests. This does not, however, exploit the stochastic nature of the regressors to the fullest extent. In the present study, it is demonstrated that the incorporation of rank analysis of covariance tools in this more complex setup (\ref{model}) yields rank tests which have better performance characteristics. To emphasize this enhanced efficiency, extensive numerical studies on simulated as as well as a real data set are carried out.
   Section 2 is devoted to the preliminary notions and description of the methods.
Section 3 deals with the partially linear model with i.i.d. nuisance 
covariates. Section 4 is devoted to rank analysis of covariance in partially linear 
models, and Sections 5 and 6 provide numerical
illustrations, both on simulated and real data.

\section{Preliminary notion}\label{section2} 
\setcounter{equation}{0}
We motivate our statistical models through an interesting case studied by Nummi and M\"{o}tt\"{o}nen (2004). They described a computer-based forest harvesting technique in Scandinavia, 
where the tree stems are converted into smaller logs and the stem height and diameter measurements are taken at fixed intervals. The harvester receives the length and diameter data at the $i$th stem point from a sensor, and a measuring and computing equipment enables a computer-based optimization of crosscutting. 
Nummi and M\"{o}tt\"{o}nen (2004) 
consider the model of regression dependence of the stem diameter measurement $y_i$ on the stem height measurement $x_i$ at the
$i$th stem point, $i=1, \ldots, n.$ 
The problem of interest is the prediction for $y_i$ and the testing of hypotheses on the parameters of the model; but both the stem diameter and the stem height contain measurement errors. On top of that the volume of the stem may not be 
linearly related to its diameter, rather it is more likely to be related to its height and the cross-section which may be roughly proportional to the square of the diameter. 

There are many other similar problems which can be described by partially linear regression models of the type (\ref{model})
where ${\mathbf x}_{i}$ is a $p$-vector covariate, $\mathbf Z_{i}$ is a $q$-vector covariate, the function $\nu(\cdot)$ is unknown, and the
model error $e_{i}$ is independent of $({\mathbf x}_{i}, \mathbf Z_{i}), \ i=1,\ldots,n.$
It means that the response variable $Y_{i}$
 depends on variable ${\mathbf x}_{i}$ in a linear
way but is still related to another independent variables $\mathbf Z_{i}$ in an unspecified form,
$i=1,\ldots,n.$ This model, along with the measurement errors model, are flexible and enable to model various situations with latent variables present.

In (\ref{model}) we assume that
the independent errors
$e_{1},\ldots,e_{n}$ are identically distributed according to an unknown
distribution function $F,$ and that
$\boldgreek\beta^{\top}=(\beta_1,\ldots,\beta_p),$
${\boldgreek\beta}^*=(\beta_0,{\boldgreek\beta}^{\top})^{\top}$ are unknown parameters. The function $\nu(\cdot)$
is unknown and $\mathbf Z_{i}$ are additional covariates; if they are unobservable, then all $\nu(\mathbf Z_{i}), \; i=1,\ldots,n$ are latent random variables. The rank tests of for this situation with unobservable $\mathbf Z_i$ are studied in Section 3. If $\mathbf Z_{i}$'s are observable, we can use this additional information even if $\nu(\cdot)$ remains unknown, and apply the methods of the rank analysis of covariance; very important is that this method is successful even when $\mathbf Z_{i}$ itself is affected by a measurement error (Section 4). 

Our interest is to find how the rank tests of hypothesis ${\mathbf H}_0$ in (\ref{H0})
behave in the described situations and to demonstrate their superiority to other methods. They are distribution free and avoid an estimation of nuisance $\nu(\cdot),$ which would always worsen the rate of convergence
of the whole procedure. The numerical study in Section 5 illustrates the good behavior of the rank tests in situations with various uncertaintes.

\section{Partially linear model with i.i.d. latent variables} 
\setcounter{equation}{0}
Consider the partially linear model (\ref{model}) and the problem of testing the hypothesis ${\mathbf H}_0: \ {\boldgreek\beta}={\mathbf 0},$ with $\beta_0$ and function $\nu(\cdot)$ unknown, the $\mathbf Z_i$ (scalar or vector random variables) unobservable. The model can be rewritten as
\begin{equation} \label{model2}
 Y_i=\beta_0+{\bf x}_{i}^{\top}{\boldgreek\beta}+e_i^*, \; e_i^*=e_i+\nu(\mathbf Z_i), \;
 i=1,\ldots,n.
\end{equation} 
The regression matrix $\mathbf{X}=\mathbf{X}_n$ in model (\ref{model}) is of order $n\times p$
with the rows $\mathbf{x}_{i}, \ i=1,\ldots,n.$ Denote $\mathbf{X}_n^0$ the matrix with the rows
$\mathbf{x}_{i}-\bar{\mathbf{x}}_n,$ $i=1,\ldots,n,$ and assume that it satisfies
\begin{eqnarray}\label{25}
&&\mathbf{Q}_n=\frac 1n\mathbf{X}_n^{0\top}\mathbf{X}_n^0=\frac 1n\sum_{i=1}^n(\mathbf{x}_{i}-\bar{\mathbf{x}}_n)
(\mathbf{x}_{i}-\bar{\mathbf{x}}_n)^{\top}\rightarrow
\mathbf{Q} \ \mbox{ as } \ \ny,\\[2mm] &&n^{-1}\max_{1\leq i\leq
n}\left\{(\mathbf{x}_{i}-\bar{\mathbf{x}}_n)^{\top}\mathbf{Q}_n^{-1}
(\mathbf{x}_{i}-\bar{\mathbf{x}}_n)\right\}\rightarrow 0 \ \mbox{
as } \ \ny\nonumber
\end{eqnarray}
where $\mathbf{Q}$ is a positive definite $p\times p$ matrix.

Assume that the distribution function $F$ of the errors $e_i$ has an absolutely continuous
density $f$ and finite Fisher information ${\mathcal I}(f)=\int_{\mathbb{R}}\left(\frac{f^{\prime}(z)}{f(z)}\right)^2dF(z)<\infty.$
Assume that $\mathbf Z_1,\ldots,\mathbf Z_n$ are i.i.d.; let $G$ be the joint distribution function of $\nu(\mathbf Z_i), \ i=1,\ldots,n.$ It is unknown, we only assume that it has an absolutely continuous density $g$ and finite Fisher information $\mathcal I(g).$ Moreover, let $H$
denote the distribution function of $e_i^*, \ i=1,\ldots,n.$ Because $e_i^*$ is more dispersed than $e_i,$ then $\mathcal I(h)\leq\mathcal I(f),$ where $\mathcal I(h)$ is the Fisher information of $H,$ with the equality if $\nu(\mathbf Z_i)=0$ with probability 1 (see H\'ajek et al. (1999)).

Let $R_1,\ldots,R_n$ be the ranks of $Y_1,\ldots,Y_n.$ 
The rank tests of
$\mathbf{H}_0: \ \boldgreek{\beta}=\mathbf{0},$ both in models (\ref{model}) and (\ref{model2}), are based on the
vector of linear rank statistics $\mathbf{S}_n\in\mathbb{R}^p,$
\begin{equation}\label{26}
  \mathbf{S}_n=n^{-1/2}\sum_{i=1}^n(\mathbf{x}_{i}-\bar{\mathbf{x}}_n)a_n(R_{i})
\end{equation}
where the scores
$a_n(i)$ are generated by nondecreasing, square integrable score function
$\varphi:(0,1)\mapsto\mathbb{R}^1$ in either of the following two
ways:
\begin{eqnarray}\label{15a}
&&  a_n(i)=\mathbb{E}\varphi(U_{n:i}),\\[2mm] 
&&
a_n(i)=\varphi\left(\frac{i}{n+1}\right), \ i=1,\ldots,n,\nonumber
\end{eqnarray}
and $U_{n:1}\leq\ldots\leq U_{n:n}$ are the order statistics
corresponding to the sample of size $n$ from the $R(0,1)$
distribution.  
The test criterion for $\mathbf{H}_0$ is the quadratic form in $\mathbf{S}_n,$
\begin{equation}\label{27}
  \mathcal{T}_n^2=(A(\varphi))^{-2} \ \mathbf{S}_n^{\top}\mathbf{Q}_n^{-1}\mathbf{S}_n
\end{equation}
where
\begin{equation}\label{17a}
A^2(\varphi)=\int_0^1(\varphi(t)-\bar{\varphi})^2dt,\quad \bar{\varphi}=\int_0^1\varphi(t)dt.
\end{equation}
and because the ranks are distribution free, its asymptotic null distribution is $\chi^2$ with $p$ degrees of
freedom, and the nonlinear regressor does not cause any bias.

On the other hand, the asymptotic distributions of $\mathcal{T}_n^2$ under the local
alternative
\begin{equation}\label{28}
  \mathbf{H}_n: \
  \boldgreek{\beta}=\boldgreek{\beta}_n=n^{-1/2}\boldgreek{\beta}^*, \quad {\mathbf 0}\neq \boldgreek{\beta}^*\in
  \mathbb{R}^p \ \mbox{ fixed,}
\end{equation}
are the noncentral $\chi^2$ distributions with generally different noncentrality parameters. The relative asymptotic efficiency of the test in the presence of the nonlinear covariate with respect to that in a genuinely linear model
is given in the following theorem: 
\begin{thm}\label{Theorem1}
Let $\mathcal{T}_n^2$ be the test criterion (\ref{27}) for $\mathbf H_0$ and $\mathcal{T}_{n0}^2$ be its special case corresponding to $P(\nu(\mathbf Z)=0)=1.$ Then
\begin{description}
	\item[(i)] Under $\mathbf H_0,$ both $\mathcal{T}_n^2$ and $\mathcal{T}_{n0}^2$ have asymptotically $\chi^2$ distribution with $p$ degrees of freedom, as $\ny.$
	\item[(ii)] The asymptotic relative efficiency of $\mathcal{T}_n^2$ with respect to $\mathcal{T}_{n0}^2$ under the local alternative (\ref{28}) is 
\begin{equation}\label{ARE}
{\rm e}(\mathcal T_n^2,\mathcal T_{n0}^2)=\left(\frac{\gamma(\varphi,h)}{\gamma(\varphi,f)}\right)^2=\left(\frac{\int_0^1h(H^{-1}(t))d\varphi(t)}{\int_0^1f(F^{-1}(t))d\varphi(t)}\right)^2
\end{equation}
where $f,F$ are the density and distribution function of $e_1$ in model (\ref{model}), $h,H$ are the same for $e_1^*$ in model (\ref{model2}), and 
where
\begin{equation}\label{17b}
\gamma(\varphi,h)=\int_0^1\varphi(t)\varphi(t,h)dt, \quad
\varphi(t,h)=-\frac{h^{\prime}(H^{-1}(t))}{h(H^{-1}(t))},
\end{equation}
and similarly for $\gamma(\varphi, f).$
\end{description} 
\end{thm}  
\textbf{Proof. } 
By H\'ajek et al. (1999), Sections V.1.5 and V.1.6, we have under ${\mathbf H}_0$ as well as under ${\mathbf H}_n$ 
\begin{equation}\label{V.1.5}
\|\mathbf Q_n^{-1/2}[{\mathbf S}_n-\widetilde{\mathbf L}_n]\|=o_p(1) \; \mbox{ as } \; \ny
\end{equation}
 where 
\begin{eqnarray*}
&&\widetilde{\mathbf L}_n=n^{-1/2}\sum_{i=1}^n({\mathbf x}_{i}-\bar{\mathbf x}_n)\varphi(H(Y_{i}))
\end{eqnarray*}
here $U_{n1},\ldots,U_{nn}$ 
are the random samples from the uniform $(0,1)$ distribution. Hence, both $\mathcal T_n^2$ and $\mathcal T_{n0}^2$
are asymptotically $\chi^2$ distributed with $p$ degrees of freedom under $\mathbf H_0.$ Under $\mathbf H_n,$ the 
asymptotic distribution of $\mathcal T_n^2$ is the noncentral $\chi^2$ with $p$ degrees of freedom and with the noncentrality parameter
\begin{equation}\label{29}
  \Delta_H=\boldgreek{\beta}^{*\top}\mathbf{Q}\boldgreek{\beta}^*  \
\frac{\gamma^2(\varphi,H)}{A^2(\varphi)},
\end{equation}
while $H\equiv F$ if $\nu(Z)=0$ with probability 1. This yields (\ref{ARE}) as the relative asymptotic efficiency (ARE) of the test $\mathcal T_n^2$ with respect to the test $\mathcal T_{n0}^2.$\hfill $\Box$\\
%

For the special case of Wilcoxon scores, it follows that ${\rm e}(\mathcal T_n^2,\mathcal T_{n0}^2)\leq 1,$ 
with the equality sign holding when $\nu(\mathbf Z) = 0$ with probability 1. Similar inequality holds for
the median test, if $f$ and $g$ [density of $\nu(\mathbf Z)$] are symmetric around 0 and $f$ is unimodal, because 
then $\gamma(\varphi,h)=h(0) \leq f(0)=\gamma(\varphi,f),$ with the equality sign 
holding when $\nu(\mathbf Z) = 0$ with probability 1. For general scores, under star-shaped ordering 
of $f$ and $h$ (Doksum (1969), Bickel and Lehmann (1979)), it follows that
${\rm e}(\mathcal T_n^2,\mathcal T_{n0}^2)\leq 1.$
If the test with score function $\varphi$ is asymptotically optimal for $f,$ i.e. if
$\varphi(t)=\varphi(t,f), \ 0<t<1,$ then ${\rm e}(\mathcal T_n^2,\mathcal T_{n0}^2)\leq \frac{{\mathcal I}(h)}{{\mathcal I}(f)}\leq 1.$ In the general case,
$${\rm e}(\mathcal T_n^2,\mathcal T_{n0}^2)\leq\frac{{\mathcal I}(h)A^2(\varphi)}{\gamma^2(\varphi,f)}.$$ 

  It may be of interest whether there is a positive lower bound to (\ref{ARE}). However, allowing 
the dispersion of $\nu(\mathbf Z)$ to be large compared to that of $e,$ it can be 
shown that under the same conditions as in above, (\ref{ARE}) can be made 
arbitrarily close to 0. Thus, too much of latent effects can affect the
efficacy of rank tests; it is similar in the parametric case if
$\sigma_{\nu(\mathbf Z)}^2 /\sigma_e^2$ is large; then the latent-effects model lose the
efficacy.

Besides the presence of a nonlinear nuisance regressor, the $Y_i$ can be further affected by an additive measurement error.   Hence, instead of $Y_i$ we observe $\widetilde{W}_i=Y_i+V_i, \; i=1,\ldots,n,$ where the random errors $V_1,\ldots,V_n$ are assumed to be i.i.d. and independent of $Y_i,~ \mathbf x_i,~ \mathbf Z_i, \; i=1,\ldots,n.$ Their distribution (say $\tilde{G}$)
is unknown, we only assume that it has an absolutely continuous density $\tilde{g}.$  Then the model (\ref{model}) can be further rewritten in the form
$$\widetilde{W}_i=\mathbf x_i^{\top}\boldgreek\beta+\tilde{e}_i, \; \tilde{e}_i=e_i+\nu(\mathbf Z_i)+V_i, \; i=1,\ldots,n.$$
Let $\tilde{R}_1,\ldots,\tilde{R}_n$ denote the ranks of $\widetilde{W}_1,\ldots,\widetilde{W}_n.$ Under $\mathbf H_0,$ they are independent and identically distributed, hence 
$$P\Big((\tilde{R}_{1},\ldots,{R}_{n})=(r_1,\ldots,r_n)\Big)=\frac{1}{n!}$$
for every permutation $(r_1,\ldots,r_n)$ of $1,\ldots,n.$ The test of $\mathbf H_0$ is then based on vector of linear rank statistics
\begin{equation}\label{linear}
\widetilde{\mathbf S}_n=n^{-1/2}\sum_{i=1}^n ({\mathbf x}_{i}-{\bar{\mathbf x}}_n)a_n(\tilde{R}_{i}).
\end{equation}
The
test criterion for $\mathbf{H}_0$ is the quadratic form in $\widetilde{\mathbf{S}}_n,$
\begin{equation}\label{27b}
  \widetilde{\mathcal{T}}_n^2=(A(\varphi))^{-2} \; (\widetilde{\mathbf{S}}_n)^{\top}\mathbf{Q}_n^{-1}\widetilde{\mathbf{S}}_n.
\end{equation}
Because $\widetilde{\mathbf S}_n$ is distribution free under ${\mathbf H}_0,$
the test based on $\widetilde{\mathcal T}_n^2$ has the same null distribution as the one based on $\mathcal T_n^2,$ and their common distribution depends on the matrix $\mathbf Q_n.$ Hence, their asymptotic  null
   distributions are the same, and as such, they have the same critical
    region, which asymptotically can be approximated by the right
   hand tail of the $\chi^2$ distribution with $p$ degrees of freedom. 
  Its asymptotic distribution under the local alternative (\ref{28}) is noncentral $\chi^2$ with $p$ degrees of freedom and the noncentrality parameter
$$\Delta_{\tilde{H}}=\boldgreek{\beta}^{*\top}\mathbf{Q}\boldgreek{\beta}^*  \;
\frac{\gamma^2(\varphi,\tilde{H})}{A^2(\varphi)}$$
where $\tilde{H}$ is the distribution function of $\tilde{e}_i=e_i+\nu(\mathbf Z_i)+V_i.$

\section{Rank analysis of covariance in partially linear models}
\setcounter{equation}{0}
Consider the model (\ref{model}) as a partially linear model with possible measurement errors. If the $Y_i$ are  observed only with measurement errors, then these errors can be absorbed in the errors $e_i$ of the model. More important is when the covariates $\mathbf Z_i$ are observed only with errors, hence we only observe $\mathbf W_i=\mathbf Z_i+\boldgreek\eta_i, \; i=1,\ldots,n.$  
Hence, model (\ref{model}) can be rewritten in the form
\begin{eqnarray}\label{4.1}
&& Y_i=\beta_0+\mathbf x_i^{\top}\boldgreek\beta+e_i^{**}\\
&& e_i^{**}=e_i+\nu(\mathbf W_i),\quad \mathbf W_i=\mathbf Z_i+\boldgreek\eta_i, \quad 1\leq i\leq n,\nonumber
\end{eqnarray}
where $Y_i, ~\mathbf x_i$ and $\mathbf W_i$ are all observable, but $\mathbf W_i$ and $e_i^{**}$ may no longer be independent. Information on this dependence is recovered through the rank analysis of covariance approach, whose invariance structure enables to prevail this dependence, and even enhaces the power of the test of $\mathbf H_0.$ 
A semiparametric approach estimating $\nu(\mathbf W)$ nonparametrically, using a suitable smoothing tool, possibly leads to a slower rate of convergence; inference on $\boldgreek\beta$ is then made in a parametric way. 

 Let $R_{ni}^{(j)}$ be the rank of $W_{ij}$ among $W_{1j},\ldots,W_{nj}, \; 1\leq i\leq n; \; 1\leq j\leq q.$ Denote $\mathbf W_i=(W_{i1},\ldots,W_{iq})^{\top}, \; 1\leq i\leq n.$ Moreover, let $R_{ni}^{(0)}$ be the rank of $Y_i$
among $Y_1,\ldots,Y_n, \; 1\leq i\leq n.$ Denote
\begin{equation}\label{4.2}
\mathbb R_n= \Big[\mathbf R_{n1},\ldots,\mathbf R_{nn}\Big]
\end{equation}
the $(q+1)\times n$ rank collection matrix, where $$\mathbf R_{ni}=\left(R_{ni}^{(0)},R_{ni}^{(1)},\ldots,R_{ni}^{(q)}\right)^{\top}, \; 1\leq i\leq n.$$
Recall that under $\mathbf H_0:~\boldgreek\beta=\mathbf 0$ are $(Y_i,\mathbf W_i)^{\top}, \; i=1,\ldots,n,$  independent identically distributed $(q+1)$-vectors,  while $Y_i$ and $\mathbf W_i$ are not necessarily independent. Denote 
$G^*(\mathbf u), \; \mathbf u\in\mathbb R^{q+1}$ the distribution function of $(Y_i,\mathbf W_i)^{\top}.$ 

Define a set of $(q+1)$ scores $a_{nj}(k), \; 1\leq k\leq n$ for $j=0,1,\ldots,q,$ in the same manner as in Section 3.
For the notational simplicity, we may take $a_{nj}(k)=a_n(k), \; 0\leq j\leq q, \; k=1,\ldots,n.$ Define the random $p$-vectors
$$\mathbf T_{nj}=\frac{1}{\sqrt{n}}\sum_{i=1}^n(\mathbf x_i-\bar{\mathbf x}_n)a_n\left(R_{ni}^{(j)}\right), \; \; 0\leq j\leq q.$$
Define the matrix $\mathbf V_n$ of order $(q+1)\times(q+1)$ with the components
$$v_{nj\ell}=\frac{1}{n-1}\sum_{i=1}^n\Big(a_n\Big(R_{ni}^{(j)}\Big)-\bar{a}_n\Big)\Big(a_n\Big(R_{ni}^{(\ell)}\Big)-\bar{a}_n\Big),\quad j,\ell=0,1,\ldots,q.$$
Under $\mathbf H_0:~\boldgreek\beta=\mathbf 0,$ the $n$ columns of $\mathbb R_n$ in (\ref{4.2}) are interchangeable with the common permutational (conditional, given the set of $n!$ possible realizations of $\mathbb R_n$) probability $\frac{1}{n!}.$ Denoting this permutation measure $\mathcal P_n,$ we have 
$$\E_{\mathcal P_n}\mathbf T_{nj}=\mathbf 0, \quad \E_{\mathcal P_n}(\mathbf T_{nj}\mathbf T_{n\ell})^{\top}=v_{nj\ell}\mathbf Q_n \; \mbox{ for } \; j,\ell=0,1,\ldots,q$$
with $\mathbf Q_n$ being the matrix defined in (\ref{25}).
Decompose the matrix $\mathbf V_n$ in the form
\begin{equation}\label{4.7}
\mathbf V_n=\left[
\begin{array}{ll}
v_{n00}~&\mathbf v_{n0}^{\top}\\[2mm]
\mathbf v_{n0}~&\mathbf V_{n11}\end{array}\right] 
\end{equation}
and put
\begin{equation}\label{4.8}
v_{n00.1}=v_{n00}-\mathbf v_{n0}^{\top}\mathbf V_{n11}^{-1}\mathbf v_{n0},
\end{equation}
$$\mathbf T_{n0:1}=\mathbf T_{n0}-(\mathbf T_n^*)^{\top}\mathbf V_{n11}^{-1}\mathbf v_{n0}$$
where
$$\mathbf T_{n}^{*}=(\mathbf T_{n1}^{\top}\ldots,\mathbf T_{nq}^{\top})^{\top}.$$
Thus, $\mathbf T_{n0:1}$ is the vector of residual rank statistics of $Y_i$'s in the regression of $\mathbf T_{n0}$
on $\mathbf T_n^*.$ Note that
\begin{eqnarray*} 
&&\E_{\mathcal P_n}\mathbf T_{n0:1}=\mathbf 0,\\[2mm]
&&\E_{\mathcal P_n}(\mathbf T_{n0:1}\mathbf T_{n0:1}^{\top})=v_{n00.1}\mathbf Q_n.
\end{eqnarray*}
This suggests the test criterion
$$\mathcal L_n^0=\frac{1}{v_{n00.1}}~\left(\mathbf T_{n0:1}^{\top}\mathbf Q_n^{-1}\mathbf T_{n0:1}\right)$$
which can be further rewritten as
$$\mathcal L_n^0=\mathcal L_n-\mathcal L_n^*$$
where 
$$
\mathcal L_n=\mathbf T_n^{\top}\mathbf Q_n^{-1}\otimes\mathbf V_n^{-1}\mathbf T_n, \quad
\mathcal L_n^*=\frac{1}{v_{n00}}~\left(\mathbf T_{n0}^{\top}\mathbf Q_n^{-1}\mathbf T_{n0}\right).
$$
Regarding the rank permutation distribution $\mathcal P_n$ described above, 
we conclude that the critical region of $\mathcal L_n^0$ can be obtained by enumerating the $n!$ possible permuted values of $\mathbb R_n$ and the corresponding values of $\mathcal L_n^0.$ Due to the permutation invariance of the pertaining components, $\mathcal L_n^0$ is \textit{permutationally distribution-free} [permutation principle of Chatterjee and Sen (1964)]. Asymp\-to\-ti\-cal\-ly, as $n\rightarrow\infty,$ the permutational distribution of $\mathcal L_n^0$ can be approximated by the $\chi^2$ distribution with $p$ degrees of freedom.

Under the local alternative (\ref{28}),
\begin{equation}\label{4.13}
\mathbf V_n\stackrel{p}{\longrightarrow}\boldgreek\Gamma=\Big[\gamma_{j\ell}\big]_{j,\ell=0}^q \quad \mbox{ as } \; \ny,
\end{equation}
the limiting rank score covariance matrix. Decompose $\boldgreek\Gamma$ analogously as in (\ref{4.7}),
\begin{equation}\label{4.7a}
\boldgreek\Gamma=\left[
\begin{array}{ll}
\gamma_{00}~&\boldgreek\gamma_{0}^{\top}\\[2mm]
\boldgreek\gamma_{0}~&\boldgreek\Gamma_{11}\end{array}\right] 
\end{equation}
and put
\begin{equation}\label{4.8a}
\boldgreek\gamma_{00.1}=\boldgreek\gamma_{00}-\boldgreek\gamma_{0}^{\top}\boldgreek\Gamma_{11}^{-1}\boldgreek\gamma_{0}.
\end{equation}
%
Note that the distribution of $\mathbf T_n^*$ does not depend on (\ref{28})  [as the $\mathbf Z_i$ are i.i.d], and hence under (\ref{28}) the shifts of $\mathbf T_{n0:1}$ and of $\mathbf T_{n0}$ coincide.
Thus, under the local alternative (\ref{28})
\begin{equation}
\label{4.15}
\mathcal L_n^0\stackrel{\mathcal D}{\longrightarrow}\chi^2_{p,\Delta_H^*},
\end{equation}
the noncentral $\chi^2$ with $p$ degrees of freedom and with noncentrality parameter 
$$\Delta_H^*=\boldgreek\beta^{*\top}\mathbf Q\boldgreek\beta^{*}~\frac{\gamma^2(\varphi,H)}{\gamma_{00.1}}$$ 
with $\gamma_{00.1}$ defined in 
(\ref{4.8a}).
This further implies
\begin{equation}\label{4.16}
\gamma_{00.1}\leq \gamma_{00}=A^2(\varphi),
\end{equation}
where the equality sign holds only when $\boldgreek\gamma_0=\mathbf 0;$ hence $\Delta_H^*$ cannot be smaller than the noncentrality parameter $\Delta_H$ (see (\ref{29})) of the analysis of variance rank test with the same score function. The asymptotic relative efficiency (ARE) of the analysis of covariance rank test relative 
to the analysis of variance rank test, based on the same score function $\varphi(.)$ is given
by 
\begin{equation}\label{ARE1}
ARE~{\rm (ANOCOVA \; vs. \; ANOVA)}=\frac{\gamma_{00}}{\gamma_{00.1}}\geq 1 \; ;
\end{equation}
hence, the analysis of covariance test is always at least as efficient as the analysis of variance test. 

Summarizing, we conclude that the standard rank tests of linear hypothesis can be used even in the presence of a nonlinear nuisance regression or if there are measurement errors in the response or in the regressor, provided all these entities are  
i.i.d. and independent of each others and of the model errors. If we use the test while ignoring these disturbances,
the probability of the error of the first kind is unchanged, while the disturbances only affect the power. If the nuisance regressors are observable, using the rank analysis of covariance still enhances the power.

\subsection{Mixed linear model}
Consider the mixed model
\begin{equation}\label{Sen2.1}
Y_{i} = \beta_0 + \mathbf x_{i}^{\top}\boldgreek\beta + \mathbf Z_{i}^{\top}\boldgreek\gamma
            + e_{i}, 
\end{equation}            
where $\mathbf Z_{i}, \; i=1,\ldots,n$ are stochastic $q$-vectors and $\boldgreek\gamma$ is an unknown parameter.
The mixed linear models with random and nonrandom covariates were studied in monographs by Khuri et al. (1998) and by Muller and Stewart (2006); 
the first one has a more theoretical flavor, while the second one focuses on detailed applications. 
However, if the random covariate is observed with an error, even the mixed linear model 
leads to the form (\ref{model}) with a nonlinear nuisance regressor.  Assume that the $\mathbf 
Z_{i}$ are not directly observable, but are subject to measurement errors
 $\boldgreek\eta_{i},$ hence the observable random vectors are
$\mathbf W_{i} = \mathbf Z_{i} + \boldgreek\eta_{i}.$ 
Assume that
the $\boldgreek\eta_{i}$ are independent of both $\mathbf Z_{i}$ and $e_{i}.$ 
Without loss of generality assume that the $\E\mathbf W_{i}=\mathbf 0.$ 
Notice that $ Y_{i}$ and 
$\mathbf W_{i}$ are independent, given $\mathbf Z_{i}.$ Hence, the conditional distribution function of $Y_i^o = Y_{i} - \beta_0 - \mathbf x_{i}^{\top}\boldgreek\beta$ given $\mathbf W_{i},$ denoted by
$f_{(y^o|\mathbf w)}(y^o | \mathbf w),$ can be written as
\begin{eqnarray}\label{Sen2.2}
 &&f_{(y^o| \mathbf w)}(y^o | \mathbf w) = \int_{\mathbb R^q} f_{(y^o, \mathbf z | \mathbf w)}(y^o, \bf z | \bf w )d\mathbf z\nonumber\\ 
&& =\int_{\mathbb R^q}f_{(y^o|\mathbf w,\mathbf z)}(y|\mathbf w,\mathbf z)f_{(\mathbf z|\mathbf w)}(\mathbf z|\mathbf w)d\mathbf z\\ 
&&=\int_{R^q} f_{(y^o| \mathbf z)}(y^o -  \boldgreek\gamma^{\top} \mathbf z ) f_{(\mathbf z |\mathbf
         w)}( \mathbf z | \mathbf w) d\mathbf z.\nonumber 
\end{eqnarray}         
If the two conditional densities are Gaussian, then
(\ref{Sen2.1}) corresponds to the linear measurement error model, with $\boldgreek\gamma$ replaced by 
$\mathbf K^{\top}\boldgreek\gamma$, where $\mathbf K$ is the
matrix of $\boldgreek\Sigma_{\mathbf z}( \boldgreek\Sigma_{\mathbf w})^{-1}.$ However, if the two densities are not Gaussian,
then (\ref{Sen2.2}) involves a nuisance 
function $\nu(\mathbf W_{i}),$ where the form of $\nu(.)$ is unspecified, depending on the unknown
densities. Note that here $\mathbf W_i$ are observable, not $\mathbf Z_i,$ hence  we will have $\nu(\mathbf W_i)$ instead of $\nu(\mathbf Z_i)$ in (\ref{model}). 
Hence note that even for the mixed linear model (\ref{Sen2.1}) if the densities are not all Gaussian, we may not have a linear model, but based on  (\ref{Sen2.2}), 
we can adapt a partially linear model as in (\ref{model}). This enables us to incorporate rank analysis of covariance tests to have better power properties. There is, however, one salient point that we need to emphasize here. Since the rank ANOCOVA test is conditionally (permutationally) distribution-free, for small to moderate sample sizes the permutation distribution needs to be enumerated to compute the permutational critical values. This task is quite manageable for small sample sizes but becomes prohibitively laborious as the sample sizes increase. Though for large samples, asymptotics work out well, for moderate to small sample sizes, to aid permutation distribution enumeration, classical resampling tools (such as jackknife or bootstrap methods) can be used. We refer to the next two sections for these refinements.

\section{Numerical illustrations}
\setcounter{equation}{0}
In order to illustrate the proposed procedures  for
finite sample situation, we have conducted a simulation study.

\noindent We considered three semiparametric partially 
linear models
\begin{equation} \label{modela}
 Y_{i}=\beta_0+{x}_{i}{\beta_1}+  w_{i} \gamma +e_{i},
\quad i=1,\ldots,n,
\end{equation}
\begin{equation} \label{modelb}
 Y_{i}=\beta_0+{x}_{i}{\beta_1}+  w_{i}^2 \delta  +e_{i},
\quad i=1,\ldots,n,
\end{equation}
\begin{equation} \label{modelc}
 Y_{i}=\beta_0+{x}_{i}{\beta_1}+  \sin(w_{i})+e_{i},
\quad i=1,\ldots,n,
\end{equation}
with $w_i=z_i+\eta_i, \quad 1\leq i\leq n,$
where $\eta_i,\quad 1\leq i\leq n$ are measurement errors.
The errors $e_i, \ i=1,\dots,n,$
were simulated  from the normal $N(0,1)$, Laplace $L(0,1)$
and Cauchy distributions, respectively.  The measurement errors $\eta_i,\quad 1\leq i\leq n,$ 
were generated  independently from the normal  $N(0,0.7)$, $N(0,2)$ and uniform $U(-1,1)$ distributions.

The design points 
$x_{1},\ldots, x_{n}$ were generated from the uniform distribution on the interval (-2,10) and
$z_{1},\ldots, z_{2}$ from the uniform distribution on the interval (-10,30). They remain fixed for all simulations under given $n$.

\noindent The following parameter values of models were used:
\begin{itemize}
\item sample sizes:  $n= 20,\,100,\,500$; 
\item $\beta_0=1$;
\item $\beta_1=-0.5,-0.4,\ldots,0,\ldots,0.4,0.5$;
\item $\gamma=3$;
\item  $\delta=-2$.
\end{itemize}

Our interest is testing the hypothesis $\mathbf{H}: \ \beta_1=0$ against alternative $\mathbf{K}: \; \beta_1\neq 0$.
We use the test criterions $\mathcal{T}_n^{2}$ in (\ref{27b}) and  $\mathcal L_n^0$ in (\ref{4.15}).
10\,000 replications of the models were simulated for each combination of the parameters and each distribution of measurement errors,
and the test criterions were then computed	for the Wilcoxon 
scores.  The level $\alpha=0.05$ test was performed every time,
the mean power of the pertaining tests was then  calculated.
Figures 1--3 compare the powers in model (\ref{modelb}) with 
standard normal distribution of errors $e_i, \; i=1, \ldots, n,$ for various sample sizes. We can see that results for small $n$, i.e. $n=20$, are not overly
good, but the results are much better for larger sample sizes.
Comparing Figures 3 and 4 shows an effect of the distribution of errors $e_i, \; i=1, \ldots, n$ in model (\ref{modelb}) with
$n=500$. 
Figures 3, 5 and 6 compare the powers for different models, i.e. for (\ref{modelb}), (\ref{modela}) and (\ref{modelc}), for sample size $n=500.$ 
Figure 7 compares the empirical powers  based on $\mathcal L_n^0$ and  $\mathcal{T}_n^{2}$ in the models (\ref{modelb}) 
and (\ref{modelc}) for single 
size $n=500$.   
Superimposing the power of the analysis of covariance rank test on the
   same for the analysis of variance rank test, we see that the analysis of covariance test performs
   better than the analysis of variance test in all cases; more prominently for large
   sample sizes and when the measurement error variance is not small
   compared to the error variance of the $e_i.$ This is perfectly in line
   with our theoretical claim in (\ref{ARE1}). When the measurement error
   variance is small, the rank covariance $\nu_{01}$ is likely to be small
   too, and hence, this supremacy of the analysis of covariance test to the analysis of variance test is
   less perceptible for n = 20 (see Fig. 1 and 2). The picture becomes
   more pronounced for larger sample sizes (Fig. 4-7). 

We  have made more extensive  simulation experiments. Particularly,
  various score functions,  design vectors,  other underlying distributions of the error terms and the measurement errors with small variance were considered. 
The results were very good for larger sample sizes, similar to Figures 3--6.  Na\-tu\-ral\-ly, the results are considerably affected by the distributions of the error terms, but on the other hand,  the  influence of  the measurement errors with small variances and of the function $\nu$ of the covariate $z$ is not so substantial.
Here, too, the analysis of covariance tests give better results.
 
\begin{figure}
\begin{center}
{\includegraphics[width=11cm]{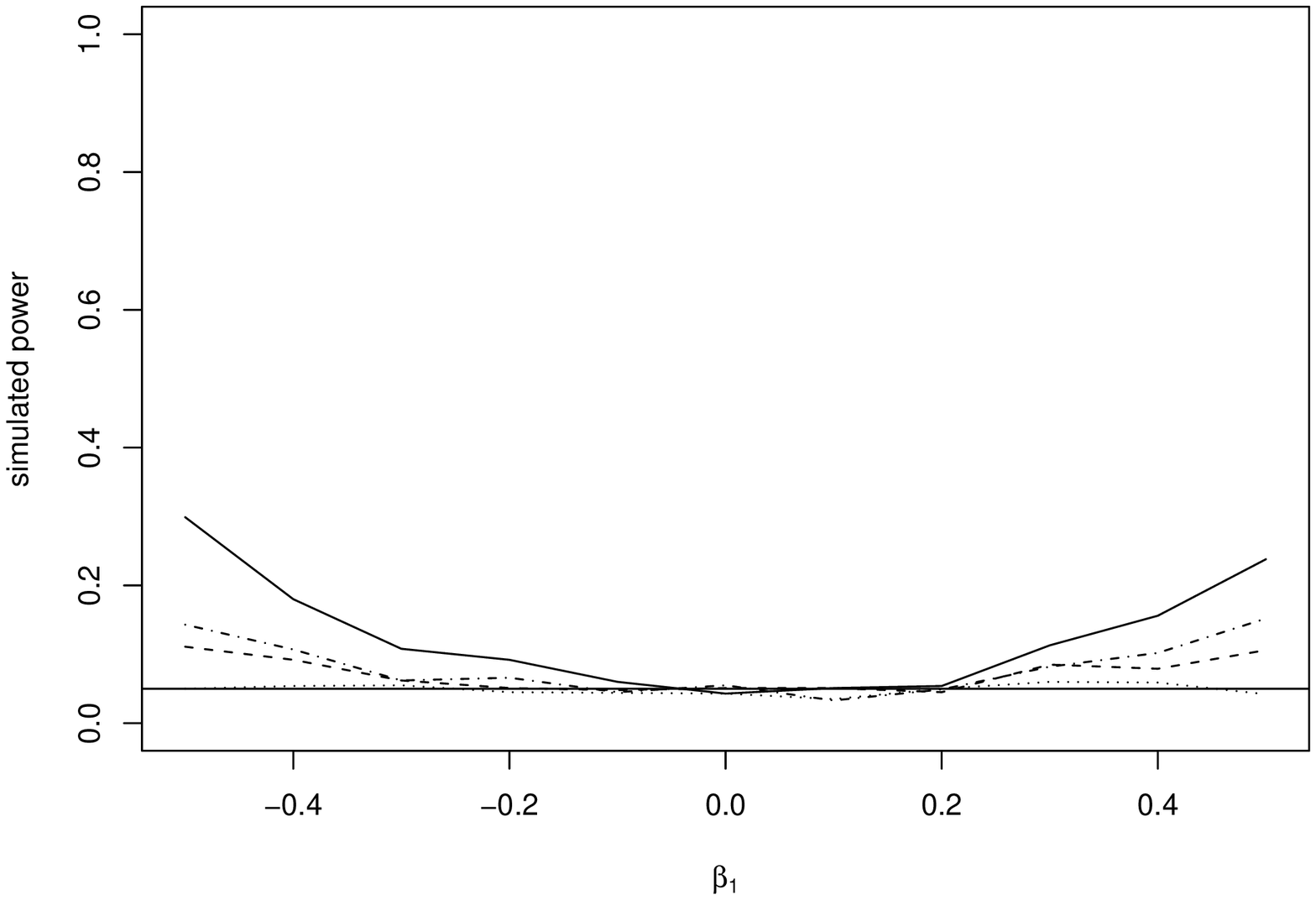}}
{\includegraphics[width=11cm]{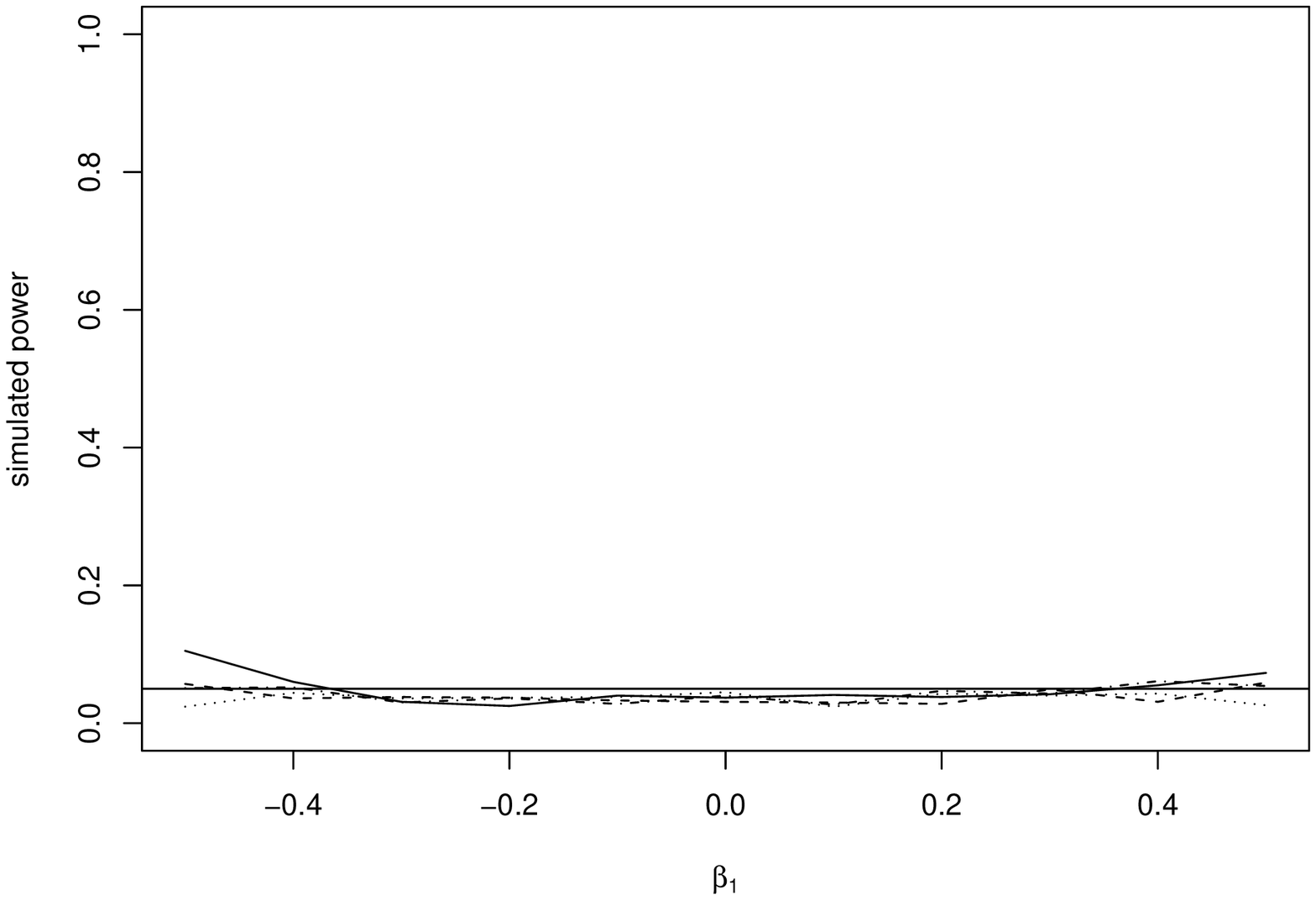}}
\caption{Empirical power of the Wilcoxon test based on $\mathcal L_n^0$(top) and  $\mathcal{T}_n^{2}$ (bottom) for $n=20$ in the model (\ref{modelb}) under the standard normal errors $e_i, \ i=1,\ldots,n$.  Solid line corresponds to the standard test, i.e. $w_i=z_i, \ i=1,\ldots,n$. The situations where $z_1,\ldots,z_n$ are affected by random errors are denoted by the dashed line (normal distribution $\mathcal N[0,0.7]$), the dotted line (normal distribution $\mathcal N[0,2]$) and dotdash line (uniform $U[-1,1]$)%
}%
\label{figure1}
\end{center}
\end{figure}

\begin{figure}
\begin{center}
{\includegraphics[width=11cm]{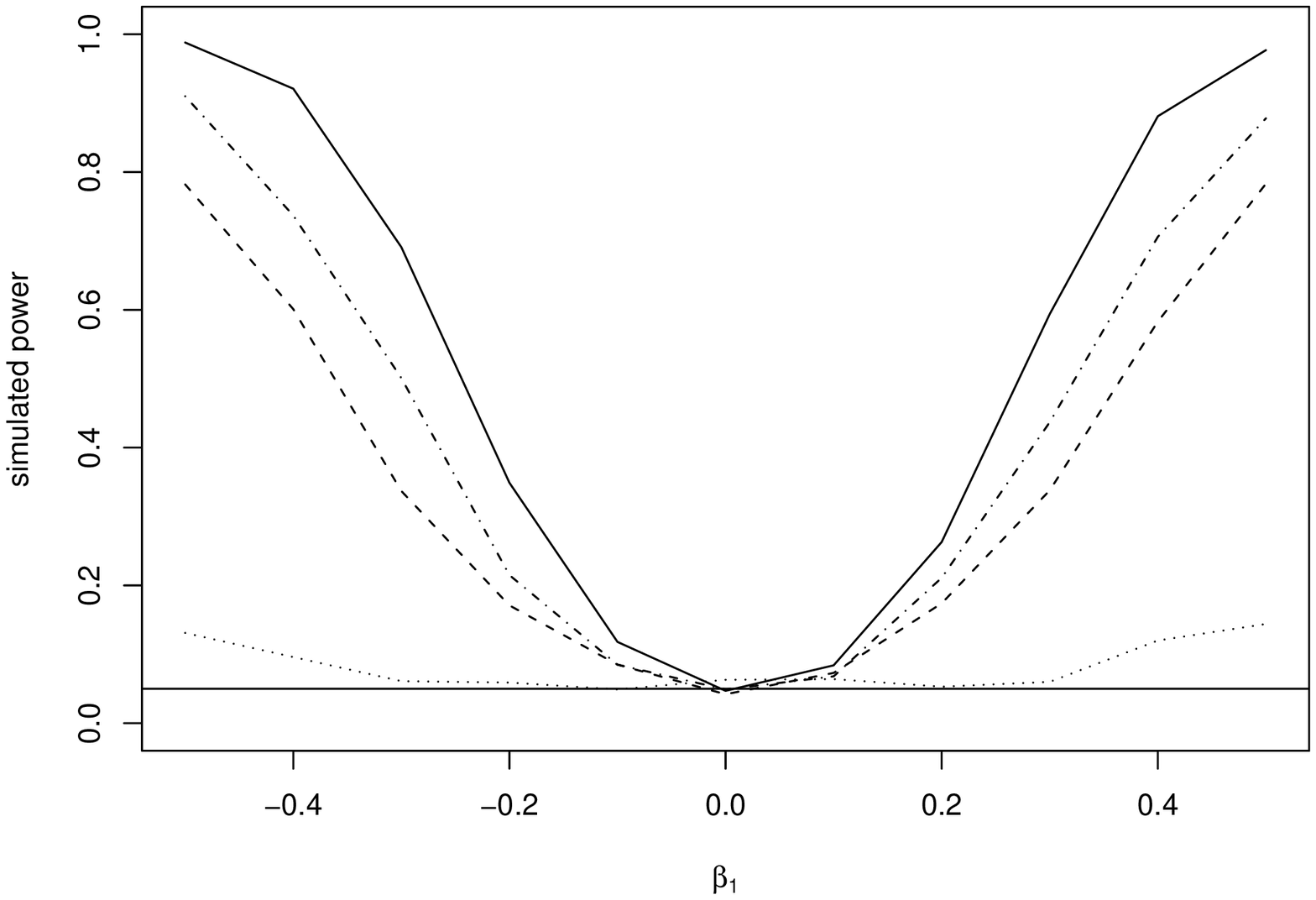}}
{\includegraphics[width=11cm]{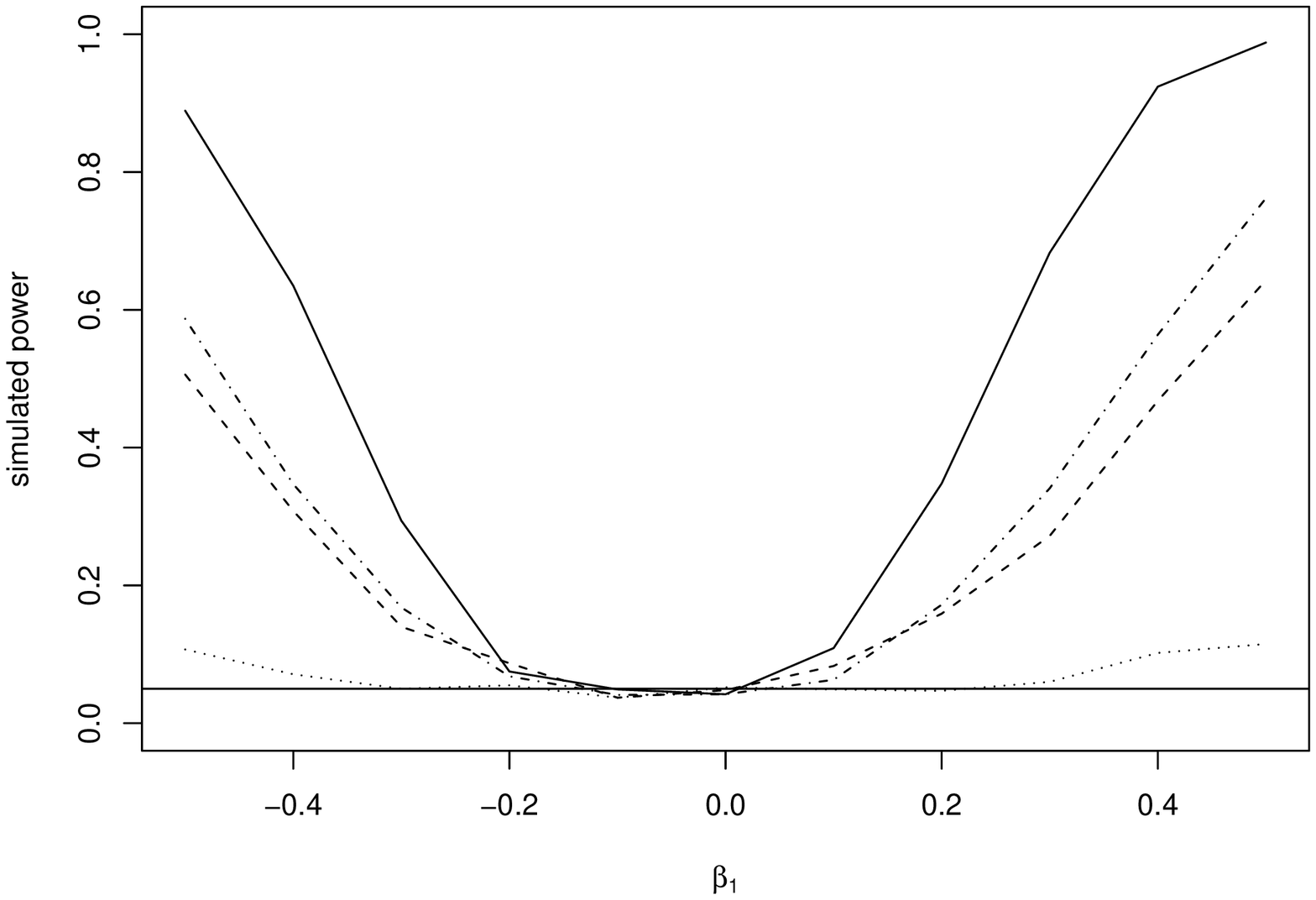}}
\caption{Empirical power of the Wilcoxon test based on $\mathcal L_n^0$(top) and  $\mathcal{T}_n^{2}$ (bottom) for $n=100$ in the model (\ref{modelb}) under the standard normal errors $e_i, \ i=1,\ldots,n$.  Solid line corresponds to the standard test, i.e. $w_i=z_i, \ i=1,\ldots,n$. The situations where $z_1,\ldots,z_n$ are affected by random errors are denoted by the dashed line (normal distribution $\mathcal N[0,0.7]$), the dotted line (normal distribution $\mathcal N[0,2]$) and dotdash line (uniform $U[-1,1]$)%
}%
\label{figure2}
\end{center}
\end{figure}

\begin{figure}
\begin{center}
{\includegraphics[width=11cm]{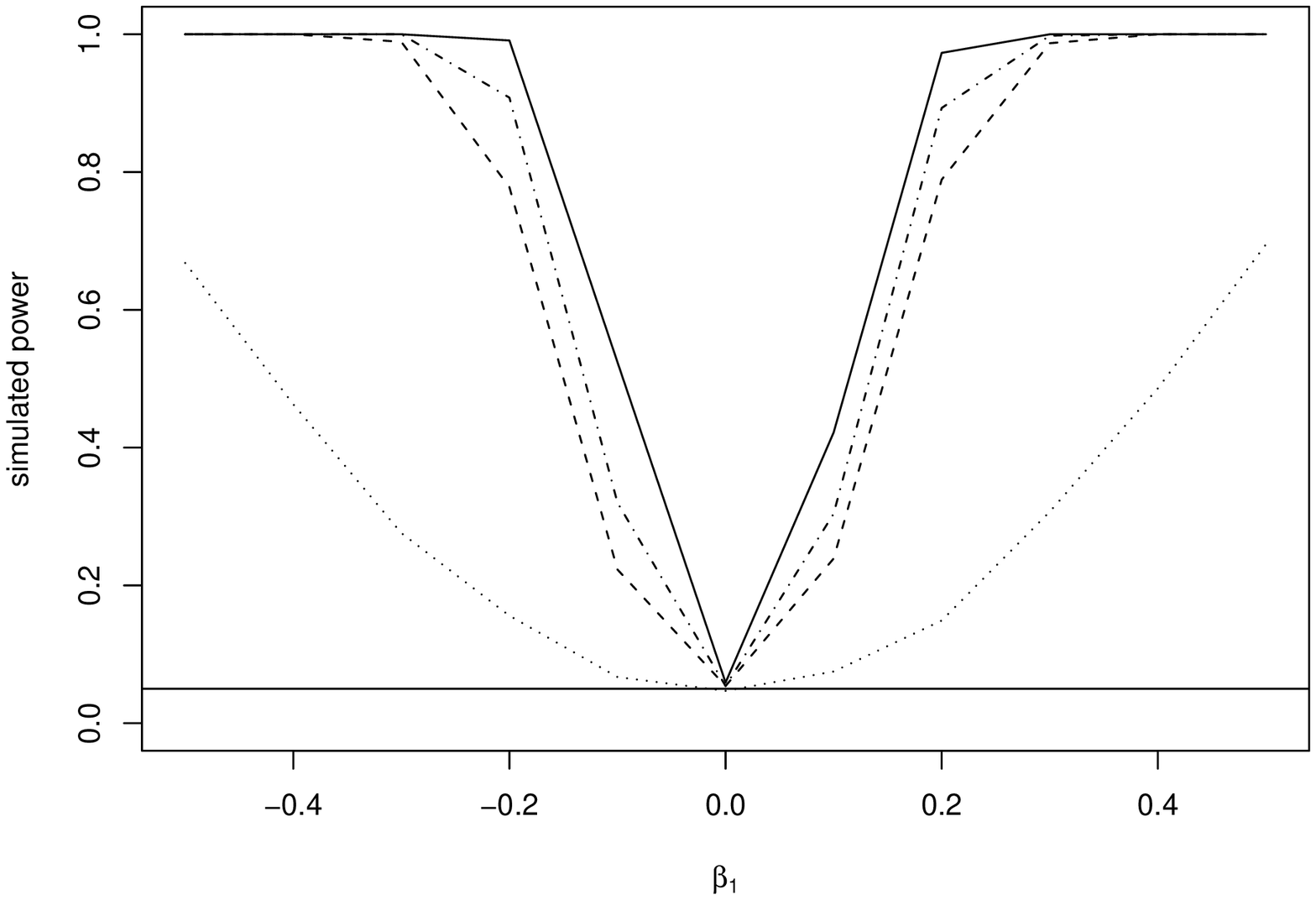}}
{\includegraphics[width=11cm]{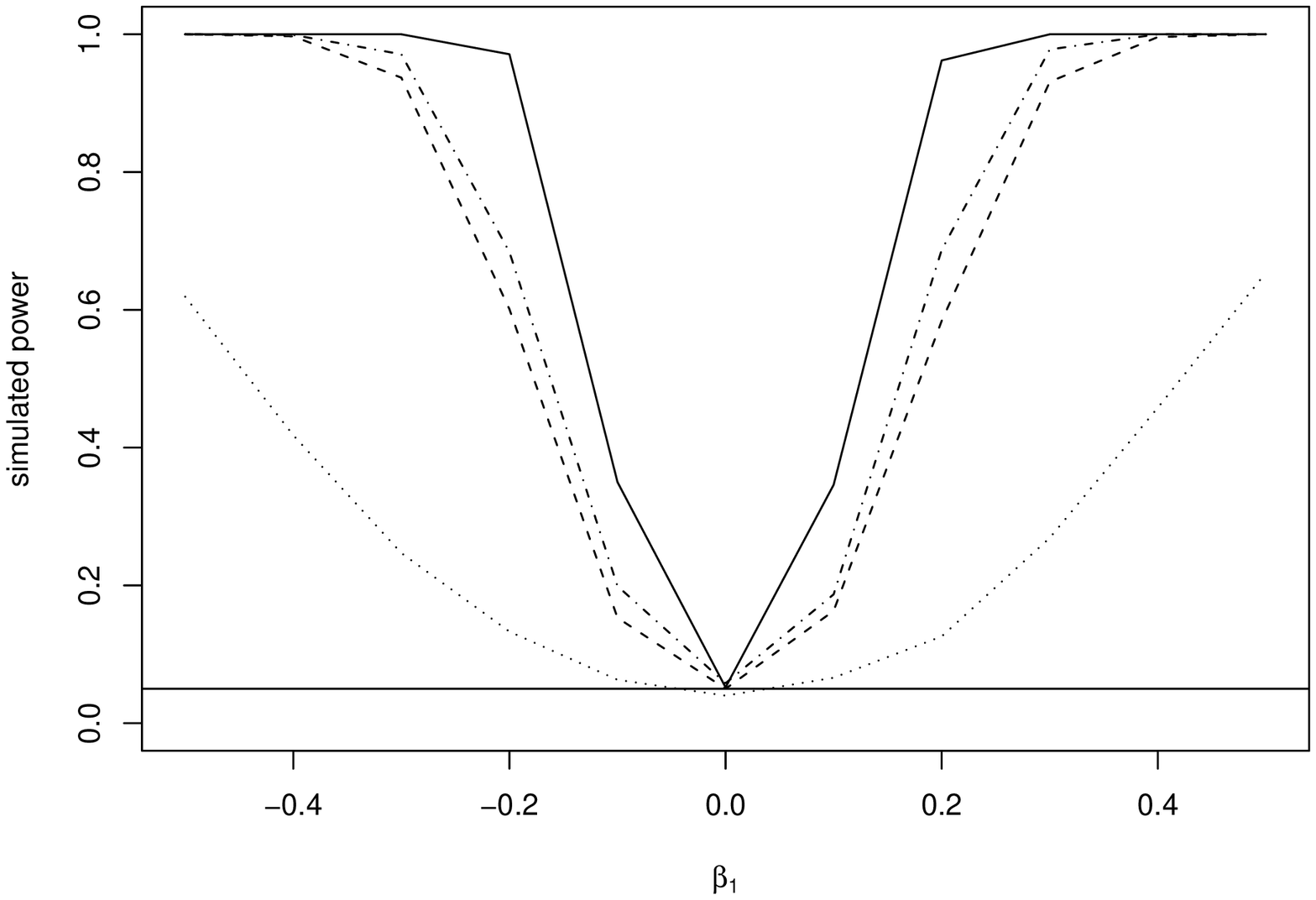}}
\caption{Empirical power of the Wilcoxon test based on $\mathcal L_n^0$(top) and  $\mathcal{T}_n^{2}$ (bottom) for $n=500$ in the model (\ref{modelb}) under the standard normal errors $e_i, \ i=1,\ldots,n$.  Solid line corresponds to the standard test, i.e. $w_i=z_i, \ i=1,\ldots,n$. The situations where $z_1,\ldots,z_n$ are affected by random errors are denoted by the dashed line (normal distribution $\mathcal N[0,0.7]$), the dotted line (normal distribution $\mathcal N[0,2]$) and dotdash line (uniform $U[-1,1]$)%
}%
\label{figure3}
\end{center}
\end{figure}

\begin{figure}
\begin{center}
{\includegraphics[width=11cm]{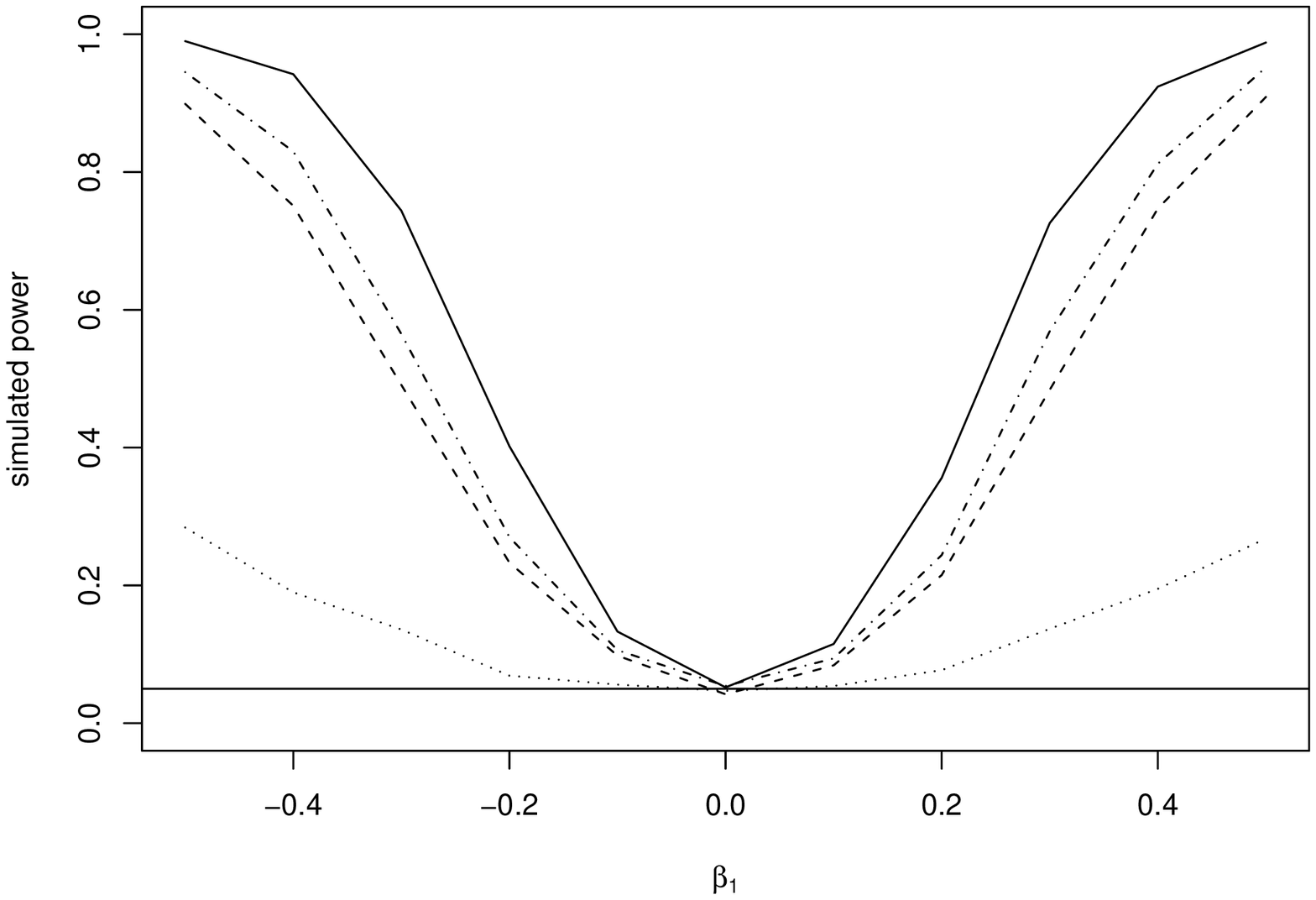}}
{\includegraphics[width=11cm]{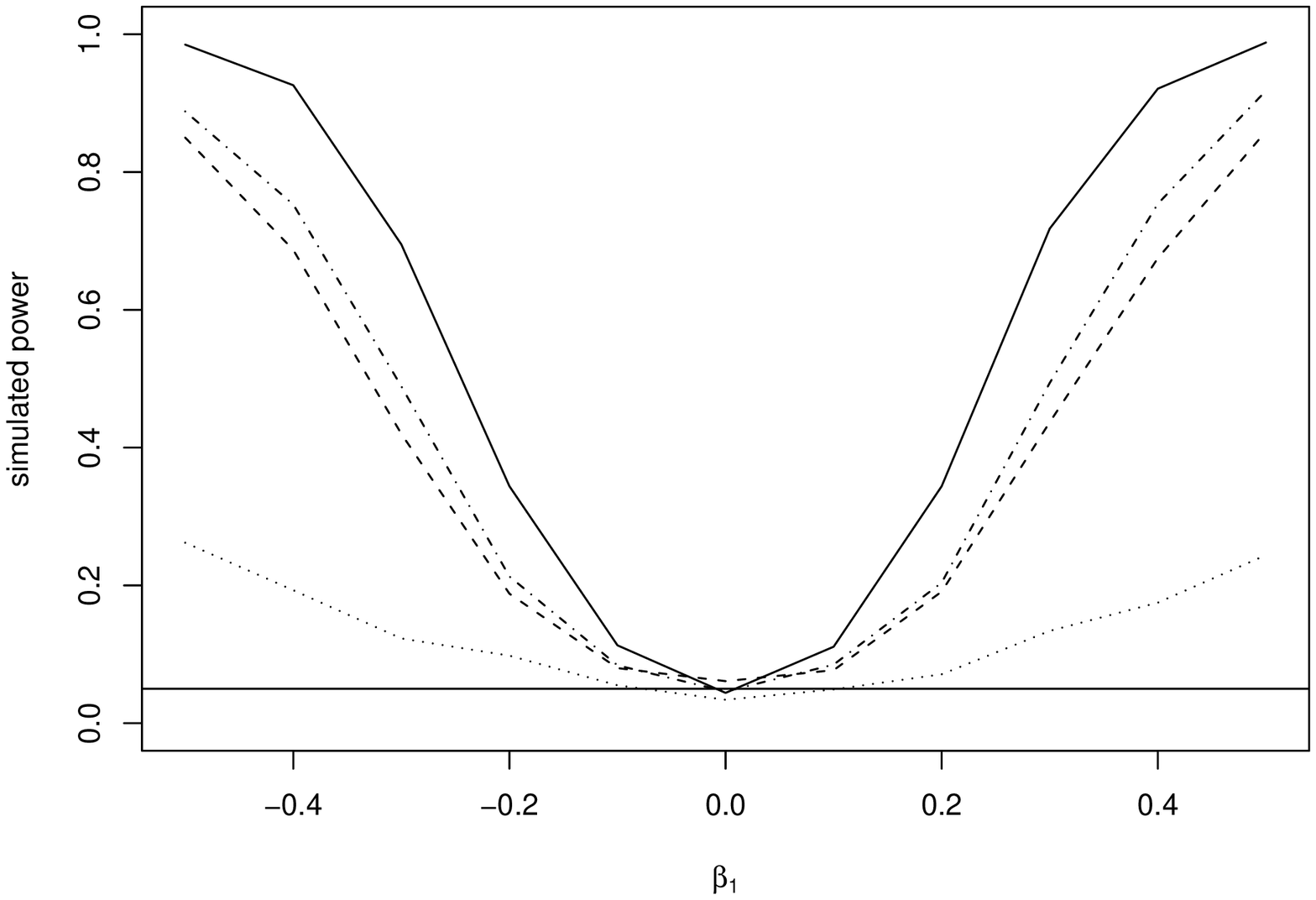}}%
\caption{Empirical power of the Wilcoxon test based on $\mathcal L_n^0$(top) and  $\mathcal{T}_n^{2}$ (bottom) for $n=500$ in the model (\ref{modelb}) under the Cauchy distributed errors $e_i, \ i=1,\ldots,n$.  Solid line corresponds to the standard test, i.e. $w_i=z_i, \ i=1,\ldots,n$. The situations where $z_1,\ldots,z_n$ are affected by random errors are denoted by the dashed line (normal distribution $\mathcal N[0,0.7]$), the dotted line (normal distribution $\mathcal N[0,2]$) and dotdash line (uniform $U[-1,1]$)%
}%
\label{figure4}
\end{center}
\end{figure}

\begin{figure}
\begin{center}
{\includegraphics[width=11cm]{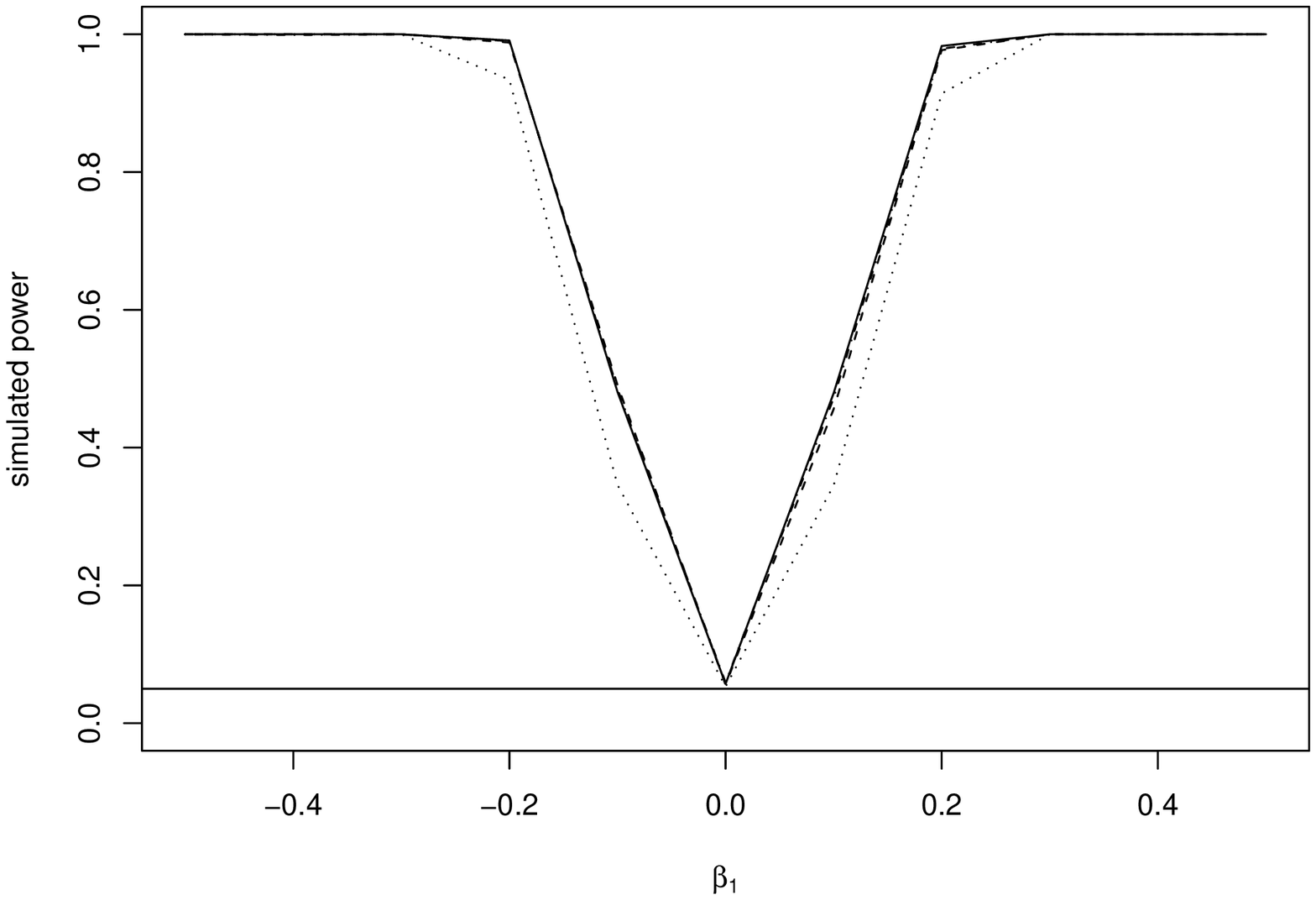}}
{\includegraphics[width=11cm]{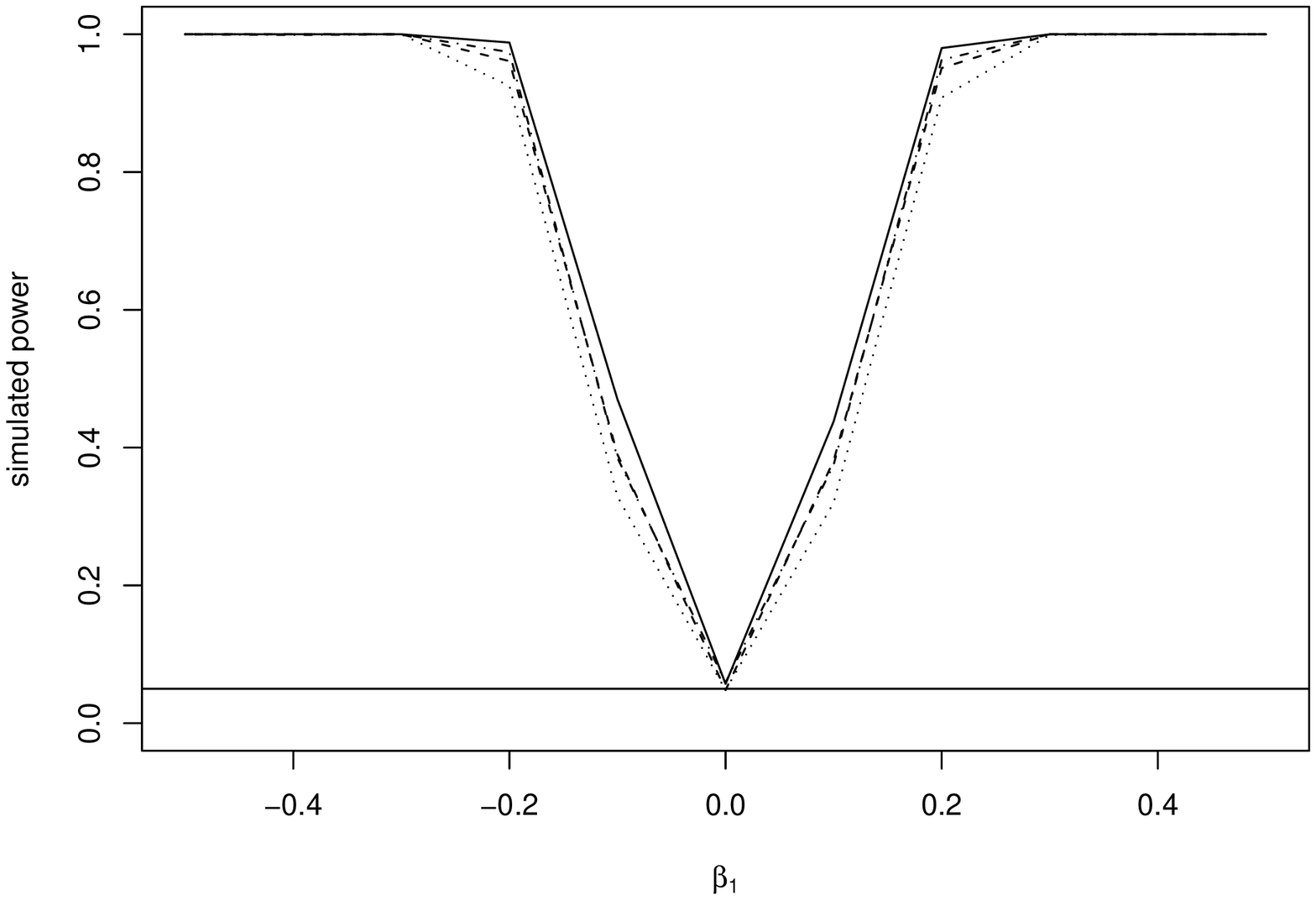}}
\caption{Empirical power of the Wilcoxon test based on $\mathcal L_n^0$(top) and  $\mathcal{T}_n^{2}$ (bottom) for $n=500$ in the model (\ref{modela}) under the standard normal errors $e_i, \ i=1,\ldots,n$.  Solid line corresponds to the standard test, i.e. $w_i=z_i, \ i=1,\ldots,n$. The situations where $z_1,\ldots,z_n$ are affected by random errors are denoted by the dashed line (normal distribution $\mathcal N[0,0.7]$), the dotted line (normal distribution $\mathcal N[0,2]$) and dotdash line (uniform $U[-1,1]$)%
}%
\label{figure5}
\end{center}
\end{figure}

\begin{figure}
\begin{center}
{\includegraphics[width=11cm]{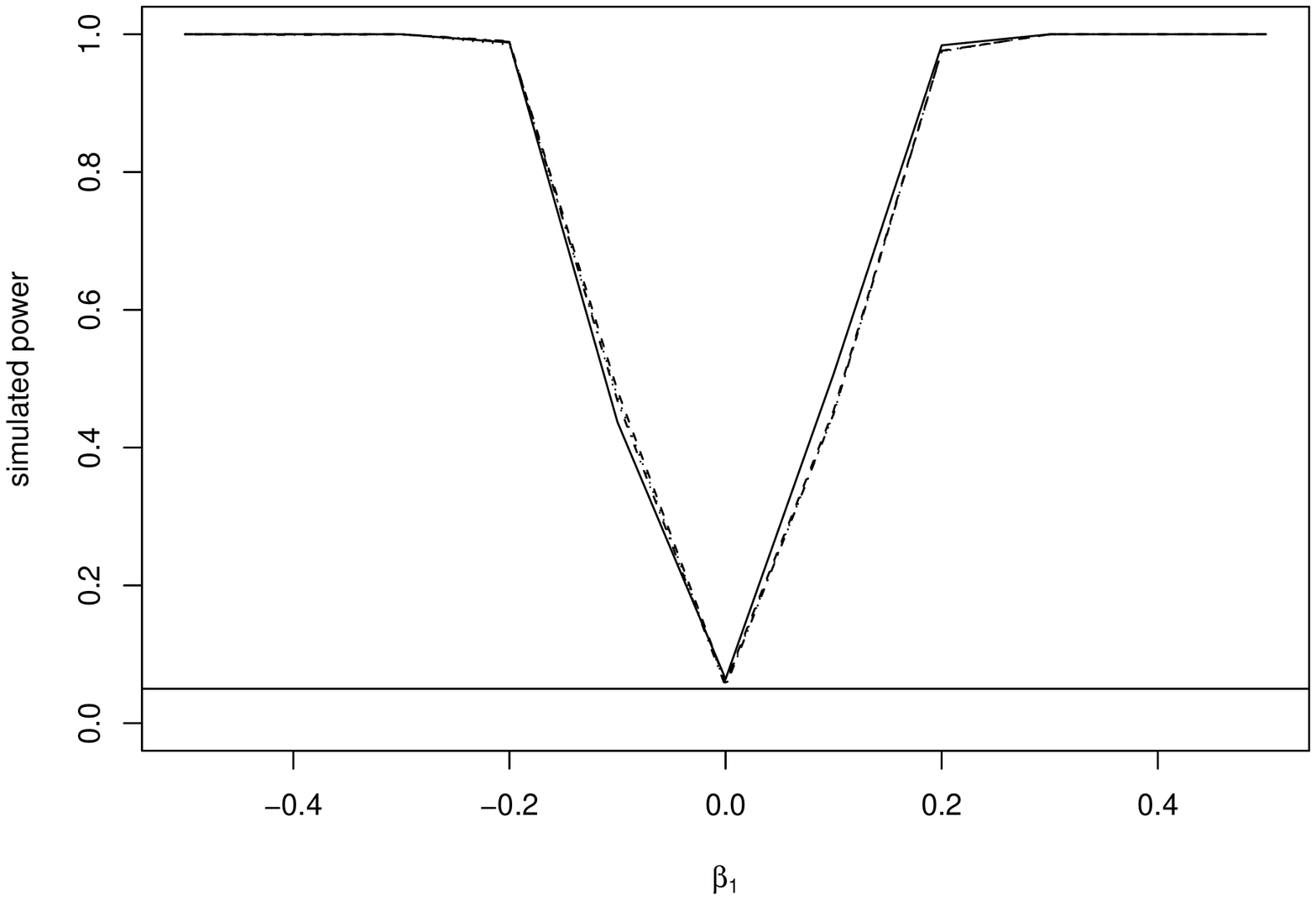}}
{\includegraphics[width=11cm]{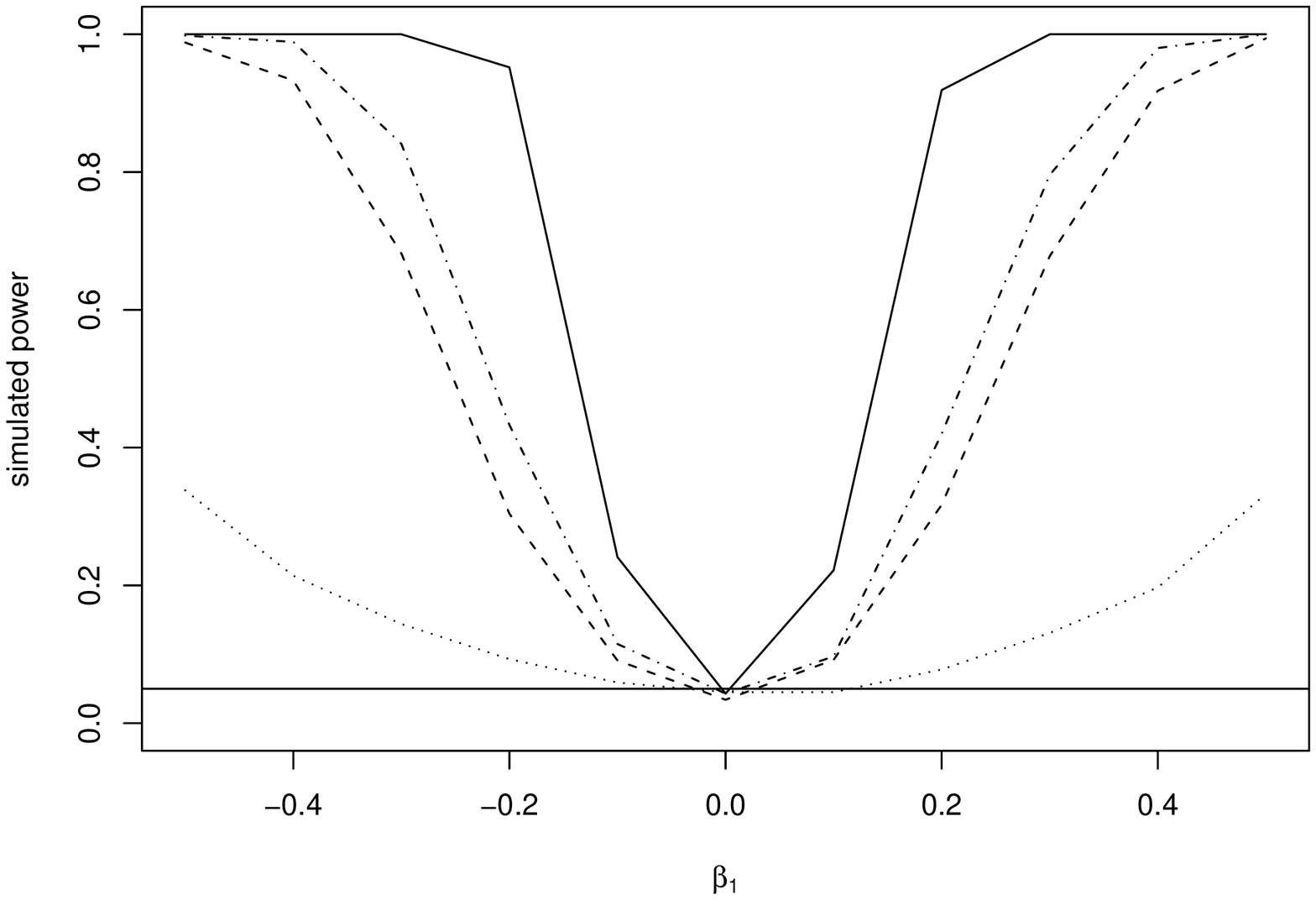}}%
\caption{Empirical power of the Wilcoxon test based on $\mathcal L_n^0$(top) and  $\mathcal{T}_n^{2}$ (bottom) for $n=500$ in the model (\ref{modelc}) under the standard normal errors $e_i, \ i=1,\ldots,n$.  Solid line corresponds to the standard test, i.e. $w_i=z_i, \ i=1,\ldots,n$. The situations where $z_1,\ldots,z_n$ are affected by random errors are denoted by the dashed line (normal distribution $\mathcal N[0,0.7]$), the dotted line (normal distribution $\mathcal N[0,2]$) and dotdash line (uniform $U[-1,1]$)%
}%
\label{figure6}
\end{center}
\end{figure}

\begin{figure}
\begin{center}
{\includegraphics[width=11cm]{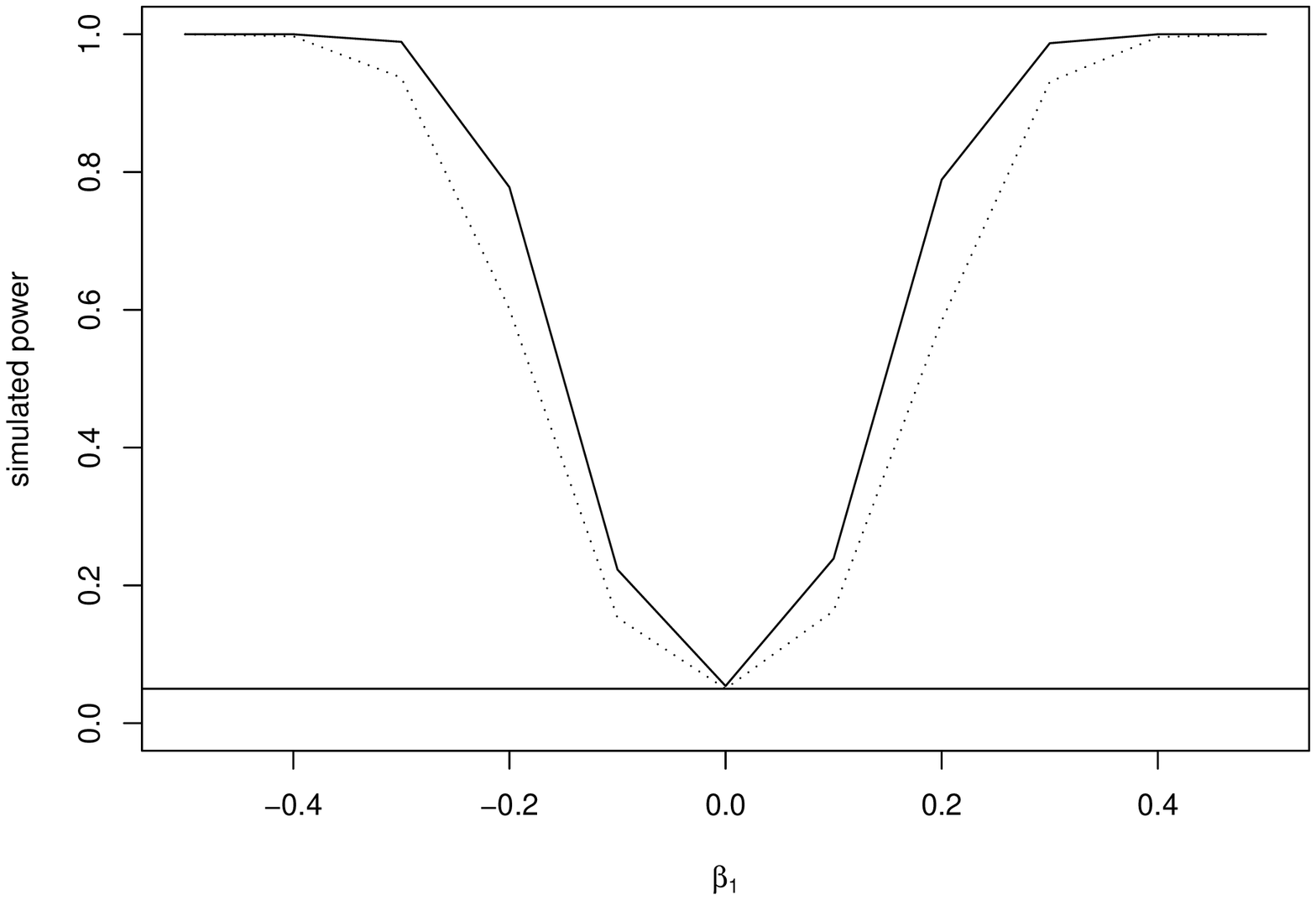}}
{\includegraphics[width=11cm]{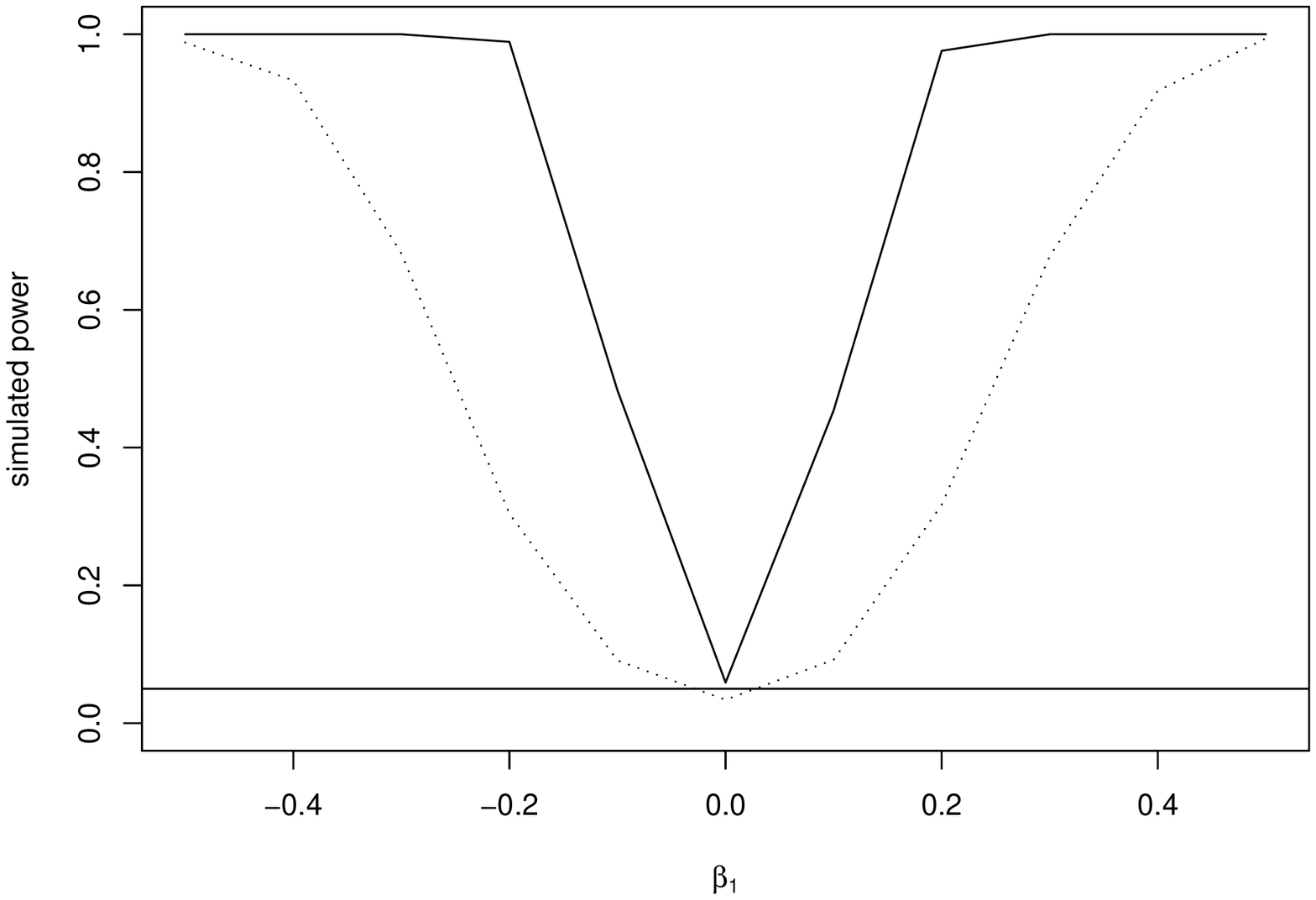}}
\caption{Comparison of empirical power  based on $\mathcal L_n^0$ (solid line) and  $\mathcal{T}_n^{2}$ (dotted line) for $n=500$ in the models (\ref{modelb}) (top)
and   (\ref{modelc}) (bottom) under the standard normal errors $e_i, \ i=1,\ldots,n.$ The covariates $z_1,\ldots,z_n$ are affected by random errors coming from normal distribution $\mathcal N[0,0.7]$%
}%
\label{figure7}
\end{center}
\end{figure}

\section{Application to the precipitation dataset}
\setcounter{equation}{0}
The test described above is applied to a  datasets of 
1-day precipitation amounts. This application makes use outputs of 
coupled atmosphere and ocean general circulation models of the NOAA Geophysical Fluid Dynamics Laboratory.
The outputs with the daily resolution are available in the form of transient climate change simulations carried out under increasing 
greenhouse gas concentrations according to prescribed emission scenarios over 1961-2100. Models 
have a horizontal resolution $2.5 \times 2.0^{\circ}$ (longitude $\times$ latitude) for South America. 

A variable of primary interest $Y$ (precipitation) is modeled using additional covariates: the time index $x$ and the southern oscillation index $Z$,
which is calculated from the monthly or seasonal fluctuations in the air pressure difference between Tahiti and Darwin.
The model under consideration has the form:
$$
 Y_{i}=\beta_0+{x}_{i}{\beta_1}+  \nu(Z_{i})+e_{i},
\quad i=1,\ldots,n,
$$
where  the $Z_{i}$ are observable but probably with measurement errors, and  $\nu(.)$ is unknown.

For each scenario gridpoint we tested the significance of  time index, i.e. $\mathbf{H}: \ \beta_1=0$ against alternative $\mathbf{K}: \; \beta_1\neq 0$. 
Table 6.1 summarizes results of testing for all 888 gridpoints and three scenarios.

\begin{center}
{\bf Table 6.1.}  {\it Rejection and non-rejection of the null hypothesis
at level $\alpha=0.05$}\\[5mm]
\begin{tabular}{|l|c|c|}
\hline
scenario&\# of rejection $\mathbf{H}$& \# of non-rejection of $\mathbf{H}$\\
\hline
scenario 1  (m21af)& 465  &  423\\  
scenario 2 (m21a2) & 576  &  312 \\
scenario 3  (m21b1  &598  &  290 \\  
\hline
\end{tabular}
\end{center} 

\section*{Acknowledgments}
The authors are grateful to the reviewer for helpful comments on the original manuscript.

\end{document}